\documentclass[11pt]{article}

\usepackage[T1]{fontenc}
\usepackage[utf8]{inputenc}
\usepackage[a4paper,margin=30mm]{geometry}
\usepackage{amsmath,amssymb,amsthm,mathtools}
\usepackage{microtype}
\usepackage[numbers,sort&compress]{natbib}
\usepackage{tikz}
\usetikzlibrary{arrows.meta}
\usepackage[hidelinks]{hyperref}
\hypersetup{pdftitle={A sharp polygonal Finsler--Hadwiger inequality},
  pdfauthor={Beniamin Bogosel}}
\numberwithin{equation}{section}
\newtheorem{theorem}{Theorem}[section]
\newtheorem{lemma}[theorem]{Lemma}
\newtheorem{corollary}[theorem]{Corollary}
\newtheorem{proposition}[theorem]{Proposition}
\theoremstyle{remark}
\newtheorem{remark}[theorem]{Remark}

\newcommand{\Dn}{\Delta_n}
\definecolor{geomblue}{HTML}{174B66}
\definecolor{geomorange}{HTML}{A04423}
\tikzset{
  diagram/.style={font=\small,line cap=round,line join=round,>=Latex},
  polygon/.style={draw=geomblue,fill=geomblue!8,line width=.8pt},
  guide/.style={draw=black!60,line width=.45pt,densely dashed},
  accent/.style={draw=geomorange,line width=1.1pt}
}

\title{A sharp polygonal Finsler--Hadwiger inequality}
\author{Beniamin Bogosel\thanks{Faculty of Exact Sciences,
  Aurel Vlaicu University of Arad, Arad, Romania.\newline
  Email: \href{mailto:beniamin.bogosel@uav.ro}{\texttt{beniamin.bogosel@uav.ro}}.}}
\date{November 24, 2026}

\begin{document}
\maketitle

\begin{abstract}
For every simple $n$-gon of area $A$ with side lengths
$a_1,\ldots,a_n$, we prove
\[
 \sum_{i=1}^n a_i^2 \ge 4\tan\frac{\pi}{n}\,A
 +c_n\sum_{i<j}(a_i-a_j)^2,
 \qquad
 c_n=1-\frac{\tan(\pi/n)}{\tan(\pi/(n-1))},
\]
with the convention $c_3=1$, so that the case $n=3$ is exactly the
Finsler--Hadwiger inequality.  The coefficient $c_n$ cannot be
improved, and equality holds only for regular polygons.  The extremal configurations are far from the regular $n$-gon: the
constant is approached by cyclic polygons in which one side shrinks to
zero while the others become the sides of a regular $(n-1)$-gon
(a doubled segment when $n=3$).  The proof combines the classical fact
that cyclic polygons maximize area for given side lengths with a
two-variable minimization using the half central angles.  We also
compute the larger optimal coefficient for nearly regular polygons,
and show that the best constant in the joint stability estimate for
side lengths and vertex radii grows like $n^2$.
\end{abstract}



\section{Introduction}

Let $P$ be a convex $n$-gon with perimeter $L$, area $A$, and side
lengths $a_1,\ldots,a_n$.  The polygonal isoperimetric inequality is
\begin{equation}\label{eq:classical}
 L^2\ge4n\tan\frac{\pi}{n}\,A,
\end{equation}
with equality exactly for the regular $n$-gon; see, for example,
\cite{Tang}.  Since $L^2\le n\sum_i a_i^2$ by the Cauchy--Schwarz
inequality, \eqref{eq:classical} implies the weaker estimate
\begin{equation}\label{eq:weitzenbock}
 \sum_{i=1}^n a_i^2\ge4\tan\frac{\pi}{n}\,A.
\end{equation}
For triangles, a classical result of Finsler and Hadwiger
\cite{FinslerHadwiger} strengthens \eqref{eq:weitzenbock} by a term
measuring how far the triangle is from equilateral:
\[
 a^2+b^2+c^2\ge4\sqrt3\,A
 +(a-b)^2+(b-c)^2+(c-a)^2.
\]
More recently, Giugiuc, Dao and Altintas
\cite[Theorem~1.1]{GiugiucDaoAltintas} found the analogue for
quadrilaterals, with the coefficient $1-1/\sqrt3$ in front of the sum
of the six squared side differences.

This led Dao \cite{DaoQuestion} to ask whether a refinement of the same
form holds for every $n$.  The coefficient proposed there,
$\tan(\pi/n)-\tan(\pi/(n+2))$, turns out to be too large already for
triangles, as Pinelis pointed out in his answer.  In this paper we
determine the optimal coefficient for every $n\ge3$
(Theorem~\ref{thm:main}).  Somewhat surprisingly, the extremal polygons
are far from regular.  The constant is approached by cyclic polygons
in which one side shrinks to zero, and the limit is a regular
$(n-1)$-gon when $n\ge4$ and a doubled segment when $n=3$.  If one
restricts to polygons close to the regular $n$-gon, the optimal
coefficient is strictly larger, namely $1/[n\cos^2(\pi/n)]$
(Proposition~\ref{prop:local}).

A second motivation comes from quantitative stability.  Indrei and
Nurbekyan \cite{IndreiNurbekyan} showed that the isoperimetric deficit
controls both the variance of the side lengths and the variance of the
distances from the vertices to their barycenter.  Indrei
\cite[Corollary~1.2]{Indrei} later added the area-weighted variance of
the barycentric angles for convex polygons.  Combined with a classical
radial-variance inequality, our sharp side-length estimate makes these
bounds fully explicit (Section~\ref{sec:stability}).
For triangles, optimal nonlinear relations between area, perimeter and
side-length dispersion were obtained in \cite{Bogosel}.

Central angles of cyclic polygons have been used before in related
problems.  Ku and Ku \cite{KuKu} derive area formulas and
isoperimetric refinements for families that include cyclic polygons
with two distinct side lengths.  Hitt and Zhang
\cite[Section~5]{HittZhang} use majorization of the central angles to
compare the perimeters of polygons inscribed in the same circle.
As far as we know, the optimal coefficient in front of
$\sum_{i<j}(a_i-a_j)^2$ was previously known only for $n=3$ and $n=4$.

\section{Main results}

A simple $n$-gon is a polygonal Jordan curve with $n$ distinct vertices,
counted with their labels even when consecutive sides are collinear.
Its area is that of the enclosed region.  A convex $n$-gon is a simple
$n$-gon that bounds a convex region.  Set
\[
 \tau_n=\tan\frac{\pi}{n}\quad(n\ge3),\qquad
 \tau_2=+\infty,
\]
and interpret $\tau_n/\tau_2$ as zero.  This convention only serves to
include the case $n=3$ in the general formula below, where it gives
$c_3=1$.

\begin{theorem}\label{thm:main}
For every convex $n$-gon $P$, $n\ge3$, with side lengths
$a_1,\ldots,a_n$ and area $A$,
\begin{equation}\label{eq:main}
 \sum_{i=1}^n a_i^2
 \ge4\tau_n A+c_n\sum_{i<j}(a_i-a_j)^2,
 \qquad
 c_n=1-\frac{\tau_n}{\tau_{n-1}}.
\end{equation}
The coefficient $c_n$ is the largest possible.  Equality holds if and
only if $P$ is regular.
\end{theorem}

Before giving the details, let us explain where $c_n$ comes from.
For fixed side lengths, the area is largest when the polygon is
inscribed in a circle.  We may therefore restrict to cyclic polygons
and describe them by their half central angles, which range over a
simplex.  If we move two of these angles while keeping their sum fixed,
the worst case is always reached either when the two angles become
equal or when one of them vanishes.  Repeating this leaves only
configurations with $m$ equal angles and $n-m$ zero ones.  The binding
case is $m=n-1$: a regular $(n-1)$-gon with one vertex counted twice,
or a doubled segment when $n=3$.

\begin{corollary}[Simple polygons]\label{cor:simple}
Every simple $n$-gon satisfies \eqref{eq:main} with the same optimal
coefficient.  Equality holds exactly for regular polygons.
\end{corollary}

\begin{proof}
The idea is to convexify $P$ without changing its side lengths.
By the Erd\H{o}s--Nagy theorem, a finite sequence of pocket flips turns
$P$ into a convex polygon $P_c$
\cite[Theorem~1.1]{DemaineGassendORourkeToussaint}.  This version of the
theorem keeps track of flat vertices in the ordered vertex list and
allows them in the convex outcome
\cite[Sections~2 and~4.2]{DemaineGassendORourkeToussaint}.
A flip reflects a boundary chain across its supporting line and leaves
the rest of the boundary in place.  Every labelled side keeps its
length, and the enclosed region gains the pocket and its mirror image
(up to boundaries), so each flip strictly increases the area.
Applying Theorem~\ref{thm:main} to $P_c$ gives \eqref{eq:main} for $P$.
If equality holds, no flip can have occurred, and we are back in the
convex equality case.  Finally, the convex sharpness family consists of
simple polygons, so the coefficient remains optimal.
\end{proof}

It is often more natural to state Theorem~\ref{thm:main} in terms of
the isoperimetric deficit and the variance of the side lengths.  Set
\begin{align*}
 D(P)&=\sum_{i<j}(a_i-a_j)^2,\qquad
 \delta_n(P)=L^2-4n\tau_n A,\\
 V(P)&=\sum_{i=1}^n\left(a_i-\frac Ln\right)^2,
 \qquad \sigma_s^2(P)=\frac1nV(P).
\end{align*}
A direct computation gives
\begin{equation}\label{eq:defectidentity}
 D(P)=n\sum_i a_i^2-L^2=nV(P),
\end{equation}
so Theorem~\ref{thm:main} can be rewritten as follows.

\begin{corollary}\label{cor:variance}
For every simple $n$-gon of positive area,
\begin{equation}\label{eq:deficit}
 \delta_n(P)\ge\gamma_nD(P),\qquad
 \gamma_n=nc_n-1=(n-1)-n\frac{\tau_n}{\tau_{n-1}}>0.
\end{equation}
Equivalently,
\begin{equation}\label{eq:variance}
 V(P)\le\frac{\delta_n(P)}{n\gamma_n},\qquad
 \sigma_s^2(P)\le\frac{\delta_n(P)}{n^2\gamma_n}.
\end{equation}
The coefficient $\gamma_n$ is maximal and the two reciprocal
coefficients in \eqref{eq:variance} are minimal.  Equality holds
exactly for regular polygons.
\end{corollary}

Indeed, subtracting the right-hand side of \eqref{eq:main} from its
left-hand side and multiplying by $n$ gives
$\delta_n(P)-(nc_n-1)D(P)$, proving equivalence.  For $n=3$,
$\gamma_3=2$.  For $n\ge4$, positivity follows because
$x\tan(\pi/x)$ is strictly decreasing for $x>2$: its derivative is
$\tan u-u\sec^2u<0$, where $u=\pi/x\in(0,\pi/2)$.

\section{Cyclic reduction and two-variable minimization}

\begin{figure}
\centering
\begin{tikzpicture}[diagram]
  \coordinate (A) at (-2.3,-.55);
  \coordinate (B) at (-1,-.55);
  \coordinate (C) at ([shift={(55:1.3)}]B);
  \coordinate (D) at ([shift={(55:1.3)}]A);
  \draw[polygon] (A)--(B)--(C)--(D)--cycle;
  \draw[->] (-.11,0)--(.47,0);
  \draw[black!32] (1.48,0) circle[radius=.91924];
  \draw[polygon] (.83,-.65) rectangle (2.13,.65);
  \node at (-.08,-1.65) {Same side lengths};

  \begin{scope}[shift={(5.0,0)}]
    \draw[black!32] (0,0) circle[radius=1.22];
    \draw[polygon] (-24:1.22)--(47:1.22)--(131:1.22)--(215:1.22)
      --(270:1.22)--cycle;
    \fill[geomorange!24] (0,0)--(-24:1.22)--(47:1.22)--cycle;
    \draw[accent] (-24:1.22)--(47:1.22);
    \draw[guide] (-24:1.22)--(0,0)--(47:1.22);
    \draw[accent,->] (-24:.38) arc[start angle=-24,end angle=47,radius=.38];
    \node at (11:.79) {$2x_i$};
    \node[geomorange] at (1.28,.29) {$a_i$};
    \fill (0,0) circle[radius=1.05pt] node[left=2pt] {$O$};
    \node at (0,-1.65) {Positive triangle};
  \end{scope}

  \begin{scope}[shift={(9.55,0)}]
    \draw[black!32] (0,0) circle[radius=1.22];
    \draw[polygon] (-56:1.22)--(-27:1.22)--(9:1.22)--(40:1.22)
      --(64:1.22)--cycle;
    \fill[geomorange!24] (0,0)--(64:1.22)--(-56:1.22)--cycle;
    \draw[accent] (64:1.22)--(-56:1.22);
    \draw[guide] (64:1.22)--(0,0)--(-56:1.22);
    \draw[accent,->] (64:.57) arc[start angle=64,end angle=304,radius=.57];
    \node at (-1.27,0) {$2x_i>\pi$};
    \fill (0,0) circle[radius=1.05pt] node[below right=1pt] {$O$};
    \node at (0,-1.65) {Negative triangle};
  \end{scope}
\end{tikzpicture}
\caption{Cyclic reduction and signed area. At left, a square maximizes
area for the rhombus's side lengths. The highlighted triangle contributes
positively to \eqref{eq:area} in the middle panel and negatively in the
right panel, where the counterclockwise gap exceeds $\pi$.}
\label{fig:cyclic}
\end{figure}
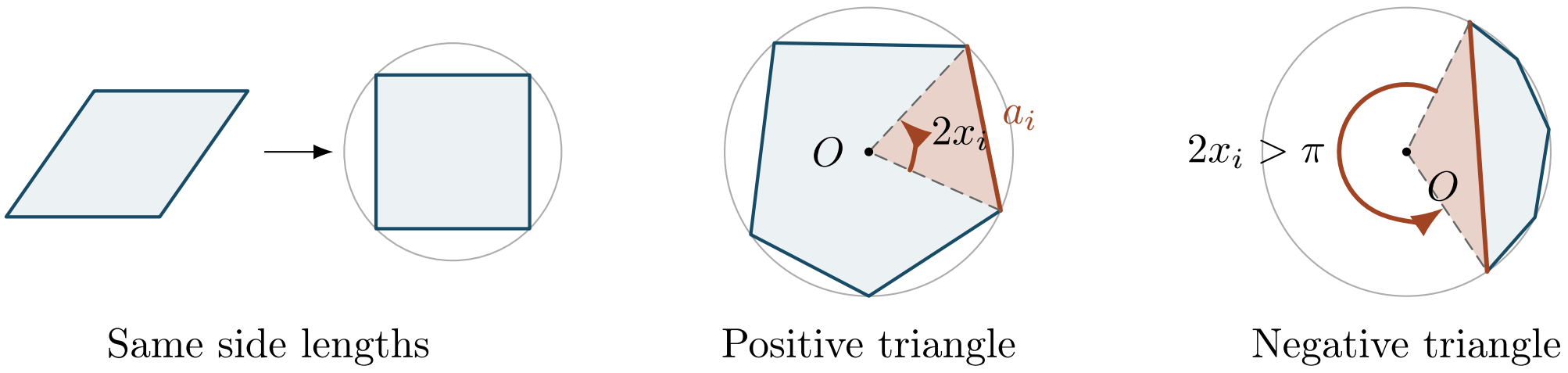

The first step is classical: among all polygons with prescribed side
lengths, the cyclic one has the largest area.  More precisely, since
$P$ has positive area, each side is strictly shorter than the sum of
the others.  Side lengths with this property, taken in boundary order,
determine a unique cyclic polygon up to isometry
\cite[Theorem~1.1]{KourimskaSkuppinSpringborn}.  By the classical
secant-polygon theorem \cite[Theorem~1.5]{KourimskaSkuppinSpringborn},
it has the largest area among polygons with these sides, and only
cyclic polygons attain it.  This remains true when $P$ has consecutive
collinear sides.  Since $A$ enters \eqref{eq:main} with a positive
coefficient, replacing $P$ by its cyclic counterpart only makes the
inequality harder, so it suffices to treat cyclic polygons
(Figure~\ref{fig:cyclic}).

Let $R$ be the circumradius, and let $2x_i$ denote the counterclockwise
angular gap from the $i$th vertex to the next.  Then
\begin{equation}\label{eq:angles}
 0<x_i<\pi,\qquad\sum_{i=1}^n x_i=\pi,
 \qquad a_i=2R\sin x_i.
\end{equation}
Note that the circumcenter need not lie inside the polygon, in which
case one of the $x_i$ exceeds $\pi/2$ (right panel of
Figure~\ref{fig:cyclic}).  The usual area formula still holds if each
triangle formed by the circumcenter and a side is counted with a sign:
\begin{equation}\label{eq:area}
 A=\frac{R^2}{2}\sum_i\sin(2x_i)
   =R^2\sum_i\sin x_i\cos x_i.
\end{equation}

Dividing the difference between the two sides of \eqref{eq:main} by
$4R^2$ and using \eqref{eq:defectidentity}, we obtain
\begin{equation}\label{eq:Phi}
 \Phi_c(x)=c\left(\sum_i\sin x_i\right)^2
 -\sum_i\left[(nc-1)\sin^2x_i+\tau_n\sin x_i\cos x_i\right]
\end{equation}
with $c=c_n$.  Theorem~\ref{thm:main} therefore amounts to showing that
$\Phi_{c_n}\ge0$.  It is convenient to allow vanishing angles, which
correspond to degenerate polygons, and to work on the compact simplex
\[
 \Dn=\left\{x\in[0,\pi]^n:\sum_i x_i=\pi\right\}.
\]

\begin{lemma}[Two-variable minimization]\label{lem:mix}
Fix $n\ge3$ and $c>0$.  There is a minimizer of $\Phi_c$ on $\Dn$
of the form
\begin{equation}\label{eq:stratum}
 x^{(m)}=\left(\underbrace{\pi/m,\ldots,\pi/m}_{m\text{ entries}},
                  0,\ldots,0\right)
 \quad\text{for some }1\le m\le n.
\end{equation}
Moreover, every minimizer with at least three positive coordinates
has all its positive coordinates equal.
\end{lemma}

\begin{figure}[tbp]
\centering
\begin{tikzpicture}[diagram]
  \foreach \shiftx/\middle/\descriptor in {
      0/30/{$z=\cos(p/2)$},
      4.2/73/{$\cos(p/2)<z<1$},
      8.4/90/{$z=1$}} {
    \begin{scope}[shift={(\shiftx,0)}]
      \coordinate (L) at (30:1.37);
      \coordinate (M) at (\middle:1.37);
      \coordinate (U) at (150:1.37);
      \draw[black!33] (10:1.37) arc[start angle=10,end angle=170,radius=1.37];
      \draw[guide] (L)--(0,0)--(U) (0,0)--(M);
      \draw[accent] (L)--(M);
      \draw[geomblue,line width=1pt] (M)--(U);
      \fill (0,0) circle[radius=1.05pt] node[below=2pt] {$O$};
      \fill[geomblue] (L) circle[radius=1.35pt]
        (U) circle[radius=1.35pt];
      \ifnum\middle=30
        \draw[accent] (M) circle[radius=2.65pt];
        \node[above right=1pt] at (M) {$x=0$};
      \else
        \fill[geomorange] (M) circle[radius=1.6pt];
      \fi
      \ifnum\middle=73
        \node at (51:.75) {$2x$};
        \node at (112:.75) {$2y$};
      \fi
      \ifnum\middle=90
        \node[above=1pt] at (M) {$x=y$};
      \fi
      \node at (0,-.71) {\descriptor};
    \end{scope}
  }
\end{tikzpicture}
\caption{A two-angle slice for consecutive gaps with the endpoints of
the arc of angle $2p$ fixed. At the two ends of the feasible interval,
the middle vertex merges with an endpoint (left) or bisects the arc (right).}
\label{fig:mixing}
\end{figure}
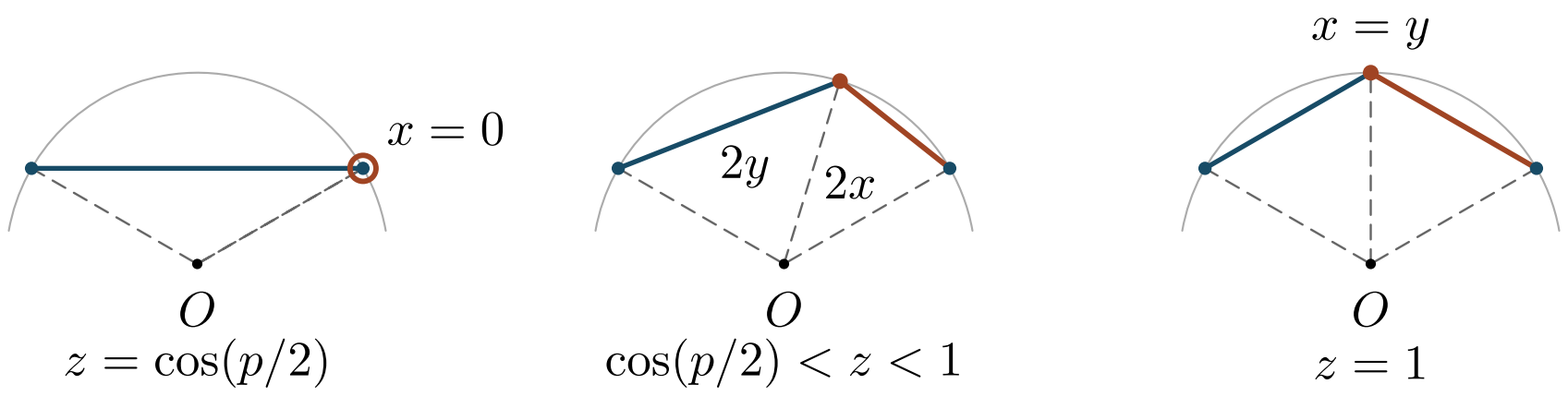

\begin{proof}
A minimum exists by compactness.  The strategy is to move two
coordinates at a time while keeping their sum fixed, and to show that
$\Phi_c$ is then concave in a suitable parameter.
Fix all but two coordinates $x,y$,
and fix their sum $p=x+y$.  If $p>0$, write
\[
 x=\frac p2+d,\qquad y=\frac p2-d,\qquad z=\cos d.
\]
Since $0\le p\le\pi$, the feasible interval is
$\cos(p/2)\le z\le1$ (Figure~\ref{fig:mixing}).  Let $r$ be the sum
of the sines of the remaining coordinates.  The pair contributions satisfy
\begin{align*}
 \sin x+\sin y&=2\sin(p/2)z,\\
 \sin^2x+\sin^2y&=1-\cos p\,(2z^2-1),\\
 \sin x\cos x+\sin y\cos y&=\sin p\,(2z^2-1).
\end{align*}
Substitution in \eqref{eq:Phi} gives
\begin{equation}\label{eq:quadratic}
 \Phi_c=\alpha z^2+\beta z+\eta,
\end{equation}
where $\eta$ is independent of $z$ and
\begin{equation}\label{eq:coefficients}
 \begin{split}
 \alpha&=4c\sin^2(p/2)+2(nc-1)\cos p-2\tau_n\sin p,\\
 \beta&=4cr\sin(p/2)\ge0.
 \end{split}
\end{equation}

The key observation is now the following.  In the variable
$t=z^2=\cos^2 d$, the slice becomes $\alpha t+\beta\sqrt t+\eta$ on
$[\cos^2(p/2),1]$, which is concave because $\beta\ge0$, and strictly
concave when $\beta>0$.  A concave function attains its minimum at an
endpoint, so we may either make one of $x,y$ vanish or make the two
equal.

Suppose a global minimizer has at least three positive coordinates.
For any two positive coordinates $x,y$, we have $p<\pi$ and $r>0$,
so $\beta>0$.  If $x\ne y$, then $t$ is interior to the interval,
which contradicts strict concavity.  Thus all positive coordinates
of such a minimizer are equal.

If a minimizer has exactly two positive coordinates, then $p=\pi$
and $r=0$.  The slice is affine on $[0,1]$, and an endpoint gives
the same minimum with one positive coordinate or two equal positive
coordinates.  A minimizer with one positive coordinate already has the
required form.  Since the coordinates sum to $\pi$, this proves
\eqref{eq:stratum}.
\end{proof}

\begin{remark}[Two positive coordinates]\label{rem:two}
On an edge of the simplex, direct substitution gives
\[
 \Phi_c(x,\pi-x,0,\ldots,0)
 =2\bigl(1-(n-2)c\bigr)\sin^2x.
\]
If $1-(n-2)c=0$, this restriction is constant.  In particular, for
$n=3$ and $c=c_3=1$, every point of such an edge has value zero,
including points whose positive coordinates are unequal.
Lemma~\ref{lem:mix} asserts only that a minimizer with equal positive
coordinates can be chosen; its characterization of every minimizer
requires at least three positive coordinates.  These edge
configurations have zero area and hence yield no additional equality
cases in Theorem~\ref{thm:main}.
\end{remark}

\section{Proof of the optimal coefficient}

By Lemma~\ref{lem:mix}, it only remains to check the sign of $\Phi_c$
at the $n$ points $x^{(m)}$.  At $x^{(1)}$, which corresponds to a
completely degenerate polygon, $\Phi_c$ vanishes.  For $m\ge2$,
evaluating \eqref{eq:Phi} at $x^{(m)}$ gives
\begin{equation}\label{eq:Phi_m}
 \Phi_c(x^{(m)})
 =m\sin^2\frac{\pi}{m}
   \left(1-c(n-m)-\tau_n\cot\frac{\pi}{m}\right).
\end{equation}
The value for $m=n$ vanishes for every $c$, as it should, since
$x^{(n)}$ is the regular polygon.  For $2\le m<n$, the value is
nonnegative exactly when
\begin{equation}\label{eq:cnm}
 c\le c_{n,m}:=\frac{1-\tau_n/\tau_m}{n-m}.
\end{equation}
The optimal coefficient is therefore $\min_{2\le m<n}c_{n,m}$.  The
next lemma shows that this minimum is attained at $m=n-1$.

\begin{lemma}\label{lem:concavity}
For $n\ge3$ and $2\le m<n$,
$c_{n,m}\ge c_{n,n-1}$, with strict inequality if $m<n-1$.
\end{lemma}

\begin{proof}
Define $h(t)=\cot(\pi/t)$ for $t\ge2$.  With $u=\pi/t$,
\begin{equation}\label{eq:concavity}
 h''(t)=-\frac{2u^3}{\pi^2\sin^2u}(1-u\cot u)<0
 \qquad(t>2),
\end{equation}
because $\tan u>u$ for $0<u<\pi/2$.  Hence $h$ is strictly
concave on $[2,\infty)$.  Since $\tau_nh(n)=1$,
\[
 c_{n,m}=\tau_n\frac{h(n)-h(m)}{n-m}
 \ge\tau_n\bigl(h(n)-h(n-1)\bigr)=c_{n,n-1}.
\]
Strict concavity gives strict inequality for $m<n-1$.
\end{proof}

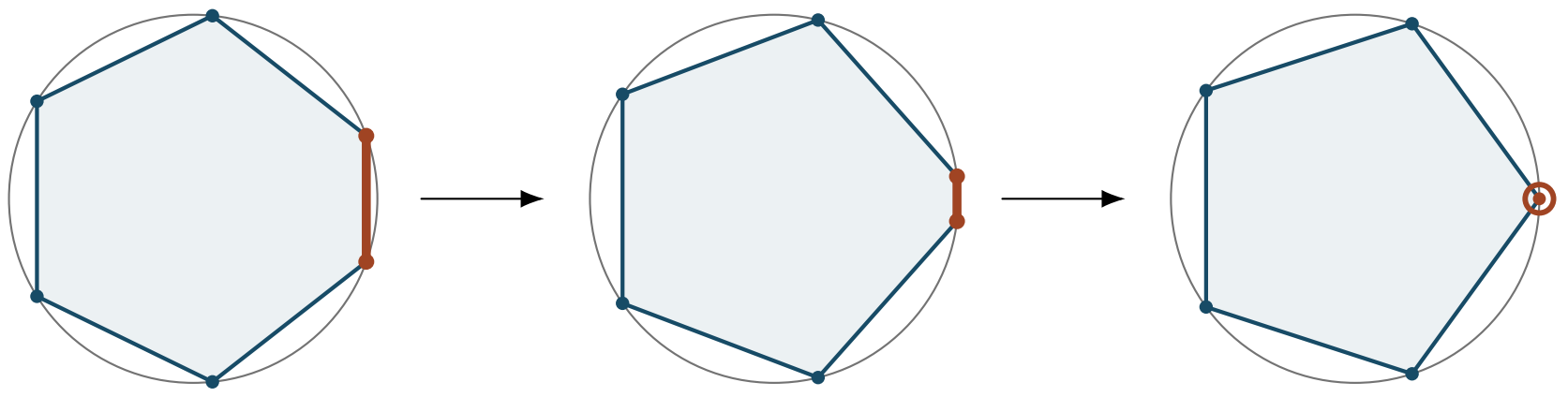
\begin{figure}[tbp]
\centering
\begin{tikzpicture}[diagram]
  \foreach \shiftx/\eps in {0/20,4.1/7} {
    \begin{scope}[shift={(\shiftx,0)}]
      \draw[black!55] (0,0) circle[radius=1.3];
      \foreach \k in {0,...,5} {
        \coordinate (v\k) at ({\eps+\k*(360-2*\eps)/5}:1.3);
      }
      \draw[polygon] (v0)--(v1)--(v2)--(v3)--(v4)--(v5)--cycle;
      \foreach \k in {1,...,4} \fill[geomblue] (v\k) circle[radius=1.35pt];
      \draw[accent,line width=1.8pt] (v5)--(v0);
      \fill[geomorange] (v5) circle[radius=1.6pt] (v0) circle[radius=1.6pt];
      \node at (0,-1.67) {$\varepsilon=\eps^\circ$};
    \end{scope}
  }
  \draw[->] (1.61,0)--(2.48,0);
  \draw[->] (5.71,0)--(6.58,0);
  \begin{scope}[shift={(8.2,0)}]
    \draw[black!55] (0,0) circle[radius=1.3];
    \draw[polygon] (0:1.3)--(72:1.3)--(144:1.3)--(216:1.3)--(288:1.3)--cycle;
    \foreach \theta in {72,144,216,288} {
      \fill[geomblue] (\theta:1.3) circle[radius=1.35pt];
    }
    \draw[accent] (0:1.3) circle[radius=2.85pt];
    \fill[geomorange] (0:1.3) circle[radius=1.3pt];
    \node at (0,-1.67) {$\varepsilon=0$};
  \end{scope}
\end{tikzpicture}
\caption{The sharpness family \eqref{eq:sharpsequence} for $n=6$.
As the thick chord of length $2\sin\varepsilon$ shrinks to zero,
the five equal chords tend to the sides of a regular pentagon.}
\label{fig:sharpness}
\end{figure}

\begin{proof}[Proof of Theorem~\ref{thm:main}]
Take $c=c_n=c_{n,n-1}>0$.  Lemma~\ref{lem:concavity} and
\eqref{eq:Phi_m} show that $\Phi_{c_n}(x^{(m)})\ge0$ for every
$m$.  Lemma~\ref{lem:mix} therefore implies
$\Phi_{c_n}\ge0$ on $\Dn$.  By the cyclic reduction at the beginning of
Section~3, this proves \eqref{eq:main} for every convex $n$-gon.

Suppose equality holds.  Equality in the area comparison first forces
$P$ to be cyclic.  All coordinates in \eqref{eq:angles} are then
positive.  Since the global minimum of $\Phi_{c_n}$ is zero,
the final assertion of Lemma~\ref{lem:mix} implies
$x_1=\cdots=x_n=\pi/n$.  Thus $P$ is regular.  Conversely, direct
substitution gives equality for regular polygons.

For sharpness, consider cyclic polygons of circumradius $1$ with
half central angles
\begin{equation}\label{eq:sharpsequence}
 x_1=\cdots=x_{n-1}=\frac{\pi-\varepsilon}{n-1},
 \qquad x_n=\varepsilon,
 \qquad 0<\varepsilon<\frac{\pi}{n}.
\end{equation}
They are strictly convex and tend to $x^{(n-1)}$ as
$\varepsilon\downarrow0$ (Figure~\ref{fig:sharpness}).
For $c>c_n$, \eqref{eq:Phi_m} gives
\[
 \Phi_c(x^{(n-1)})
 =-(n-1)\sin^2\frac{\pi}{n-1}\,(c-c_n)<0.
\]
By continuity, $\Phi_c$ is still negative for small $\varepsilon>0$, so
\eqref{eq:main} with $c$ in place of $c_n$ fails for these genuine
convex polygons.  Hence $c_n$ is optimal.
\end{proof}

Geometrically, for $n\ge4$ the sequence \eqref{eq:sharpsequence}
converges to a regular $(n-1)$-gon with one repeated vertex.  For $n=3$,
it converges to a segment traversed in both directions, with limiting
side lengths $(2,2,0)$.  In either case the limiting side-length defect
is $4(n-1)\sin^2(\pi/(n-1))>0$, so the limit is not equilateral.  The
same sequence therefore also shows that the constants in
Corollary~\ref{cor:variance} are optimal.

\section{Local sharpness and joint stability}

\subsection{Special cases and asymptotics}

For $n=3$, we have $c_3=1$ and $\gamma_3=2$, recovering the
Finsler--Hadwiger inequality.  For $n=4$,
\[
 c_4=1-\frac1{\sqrt3},\qquad
 \gamma_4=3-\frac4{\sqrt3}.
\]
To compare directly with \cite[Theorem~1.1]{GiugiucDaoAltintas},
write $Q=a^2+b^2+c^2+d^2$ and
$T=ab+ac+ad+bc+bd+cd$.  Since $D=3Q-2T$, \eqref{eq:main}
becomes
\[
 A\le\frac{T}{3+\sqrt3}-\frac{Q}{2(\sqrt3+1)^2},
\]
which is precisely their formulation.

For large $n$, the sharp coefficient behaves like $1/n$.  More
precisely, a Taylor expansion of the tangent gives
\[
 \gamma_n=\frac{2\pi^2}{3n^2}+\frac{\pi^2}{3n^3}+O(n^{-4}),
 \qquad
 c_n=\frac1n+\frac{2\pi^2}{3n^3}+O(n^{-4}).
\]
The coefficients for the two normalizations of side variance satisfy
\begin{equation}\label{eq:asymptotic}
 \frac1{n\gamma_n}\sim\frac{3n}{2\pi^2},\qquad
 \frac1{n^2\gamma_n}\longrightarrow\frac3{2\pi^2}.
\end{equation}

\subsection{The local coefficient near a regular polygon}

\begin{proposition}\label{prop:local}
As convex $n$-gons $P$ with $D(P)>0$ approach a regular $n$-gon
up to similarities, the sharp side-length quotient has lower limit
\[
 \liminf_{P\to P_{\mathrm{reg}}}
 \frac{\sum_i a_i^2-4\tau_n A}{D(P)}
 =c_n^{\mathrm{loc}}:=\frac{1}{n\cos^2(\pi/n)}>c_n.
\]
\end{proposition}

\begin{proof}
For a cyclic polygon write $x_i=\pi/n+h_i$, with $\sum_i h_i=0$.
Writing $H=\sum_i h_i^2$, Taylor expansion in \eqref{eq:Phi} gives
\[
 \Phi_c(x)=\bigl(1-nc\cos^2(\pi/n)\bigr)H+O(H^{3/2}).
\]
Indeed, the quadratic terms in $\sum_i\sin x_i$, $\sum_i\sin^2x_i$
and $\sum_i\sin x_i\cos x_i$ are respectively
$-\tfrac12\sin(\pi/n)H$, $\cos(2\pi/n)H$ and $-\sin(2\pi/n)H$.
By the definition of $\Phi_c$, $\Phi_0=(\sum_i a_i^2-4\tau_n A)/(4R^2)$
and $\Phi_0-\Phi_1=D(P)/(4R^2)$. Thus the quotient tends to
$1/[n\cos^2(\pi/n)]$ along cyclic polygons approaching the regular one.

For an arbitrary sequence approaching a regular polygon, normalize
$L=1$ so that $a_i\to 1/n$.  The equation
$\sum_i\arcsin(a_i/(2R))=\pi$ has a solution near
$R=1/[2n\sin(\pi/n)]$ for all such side lengths close enough to $1/n$.
This follows from the implicit function theorem, since every arcsine
argument is strictly between zero and one there and the derivative in
$R$ is negative.  The resulting cyclic polygon has half-angles tending
to $\pi/n$.  Its area is maximal at the given side lengths, so its quotient
is a lower bound for that of the original polygon.  Cyclic perturbations
of the regular polygon realize the limiting value.

For the strict comparison, fix $x_3=\cdots=x_n=\pi/n$ and $x_1+x_2=2\pi/n$.
At $c=c_n$, let $f(t)$ denote this slice in $t=\cos^2 d$.  The proof
of Lemma~\ref{lem:mix} gives strict concavity.  Since both endpoint
values are nonnegative and $f(1)=0$, every interior $t_0$ has
$f(t_0)>0$.  Thus the left derivative satisfies
\[
 f'_-(1)\le\frac{f(1)-f(t_0)}{1-t_0}<0.
\]
Since $t=1-d^2+O(d^4)$ and $H=2d^2$ on this slice, the quadratic
coefficient $1-nc_n\cos^2(\pi/n)=-f'_-(1)/2$ is strictly positive.
\end{proof}

The gap between the local and the global coefficients is small for
large $n$.  Expanding $\sec^2(\pi/n)$ and comparing with the expansion
of $c_n$ gives $c_n^{\mathrm{loc}}-c_n=\pi^2/(3n^3)+O(n^{-4})$.

\subsection{Recovery of side and radial stability}\label{sec:stability}

Let $v_1,\ldots,v_n$ be the vertices of $P$, let
$g=n^{-1}\sum_i v_i$ be their barycenter, and set
\[
 r_i=|v_i-g|,\qquad
 \bar r=\frac1n\sum_i r_i,\qquad
 \sigma_r^2(P)=\frac1n\sum_i(r_i-\bar r)^2.
\]
The classical radial-variance inequality of Fisher, Ruoff and
Shilleto \cite[p.~35]{FisherRuoffShilleto}, in the form recorded by
Indrei and Nurbekyan \cite[Section~3.1, equations~(3.1)--(3.2)]{IndreiNurbekyan},
gives for simple polygons
\begin{equation}\label{eq:radialclassical}
 8n^2\sin^2\frac{\pi}{n}\,\sigma_r^2(P)
 \le\delta_n(P)+n^2\sigma_s^2(P).
\end{equation}
The right-hand side of \eqref{eq:radialclassical} still contains the
side-length variance, which is exactly what
Corollary~\ref{cor:variance} controls.  Together they give a fully
explicit joint stability estimate.

\begin{corollary}\label{cor:joint}
For every simple $n$-gon,
\begin{equation}\label{eq:joint}
 \sigma_s^2(P)+\sigma_r^2(P)\le C_n\delta_n(P),
 \qquad
 C_n=\frac1{n^2\gamma_n}
       +\frac{1+\gamma_n^{-1}}{8n^2\sin^2(\pi/n)}.
\end{equation}
\end{corollary}

\begin{proof}
By \eqref{eq:variance},
$n^2\sigma_s^2(P)\le\gamma_n^{-1}\delta_n(P)$.
Substituting in \eqref{eq:radialclassical} gives
\[
 \sigma_r^2(P)\le
 \frac{1+\gamma_n^{-1}}{8n^2\sin^2(\pi/n)}\,\delta_n(P).
\]
Adding the side-length bound proves \eqref{eq:joint}.
\end{proof}

In particular, Corollary~\ref{cor:joint} provides an explicit constant
in the estimate of Indrei and Nurbekyan
\cite[Corollary~1.3]{IndreiNurbekyan}.  Using \eqref{eq:asymptotic} and
$n\sin(\pi/n)\to\pi$, we see that $C_n$ grows quadratically:
\begin{equation}\label{eq:jointasymptotic}
 C_n\sim\frac{3n^2}{16\pi^4}.
\end{equation}
The next proposition shows that this growth rate cannot be improved.
\begin{proposition}\label{prop:joint-order}
Let $C_n^\star$ be the infimum of the constants $C$ for which
$\sigma_s^2(P)+\sigma_r^2(P)\le C\delta_n(P)$ holds for all simple
$n$-gons.  For some absolute $\kappa>0$ and all
sufficiently large $n$,
\[
 \kappa n^2\le C_n^\star\le C_n=O(n^2).
\]
The lower bound holds even when the inequality is restricted to
strictly convex cyclic $n$-gons.
\end{proposition}

\begin{proof}
The upper bound follows from \eqref{eq:jointasymptotic}.  For the
lower bound, the idea is to crowd about half of the vertices near one
vertex of a regular polygon.  The resulting polygon is almost a regular
polygon with fewer sides, so its deficit is small, but the distances
from its vertices to the barycenter are very uneven.  Concretely,
fix $n\ge6$, set $m=\lfloor n/2\rfloor$, and take the
regular $m$-gon on the unit circle with one vertex at $1$. Add $n-m$
vertices on the circle near $1$. As these vertices converge to $1$,
the perimeter tends to $L_m=2m\sin(\pi/m)$ and the area tends to
$L_m^2/(4m\tau_m)$.  The limiting deficit is therefore
\[
 L_m^2\left(1-\frac{n\tau_n}{m\tau_m}\right)
 \sim\frac{4\pi^4}{n^2}.
\]
The limiting vertex barycenter is $g=(n-m)/n$. The $n-m$ repeated
vertices have radii $1-g\le1/2$, and at least $\lfloor m/2\rfloor$
regular vertices lie in the closed left half-plane and have radii
at least $1$.  The pairwise formula for the variance gives
\[
 \sigma_r^2\ge\frac{(n-m)\lfloor m/2\rfloor}{4n^2}\ge\frac1{56}.
\]
For each fixed $n$, choose the extra vertices distinct and sufficiently
close to $1$ that $\sigma_r^2\ge1/112$ and the deficit differs from
its limiting value by at most $n^{-2}$.  This gives strictly convex
cyclic polygons with $\delta_n=O(n^{-2})$ and a uniform positive lower
bound on $\sigma_r^2$, proving the required lower bound.
\end{proof}

For non-convex polygons, the deficit also controls how far $P$ is from
its convexification.  Let $P_c$ be as in the proof of
Corollary~\ref{cor:simple}.
Because the enclosed regions are nested and the side lengths are
preserved,
\[
 |P\mathbin{\triangle}P_c|=A(P_c)-A(P)
 =\frac{\delta_n(P)-\delta_n(P_c)}{4n\tau_n}
 \le\frac{\delta_n(P)}{4n\tau_n}.
\]
Combining this with \eqref{eq:joint} gives
\[
 |P\mathbin{\triangle}P_c|+\sigma_s^2(P)+\sigma_r^2(P)
 \le\left(\frac1{4n\tau_n}+C_n\right)\delta_n(P),
\]
recovering Indrei and Nurbekyan's Theorem~1.1 for simple polygons and
their Corollary~1.2 \cite{IndreiNurbekyan} with an explicit constant.

\begin{remark}[Optimal order of the deficit]\label{rem:exponent}
Along nonequilateral cyclic perturbations of a regular $n$-gon at fixed
perimeter, the computation in Proposition~\ref{prop:local} gives
\[
 \frac{\delta_n(P)}{n^2\sigma_s^2(P)}
 \longrightarrow nc_n^{\mathrm{loc}}-1=\tan^2(\pi/n)>0.
\]
Since this ratio has a positive finite limit as $\delta_n(P)\to0$, no
bound of the form $\sigma_s^2\le C\,\delta_n^{\theta}$ with $\theta>1$
can hold.  In this sense, the linear dependence on the deficit is
optimal.
\end{remark}

\section*{Statements and Declarations}

\noindent\textbf{Funding.}
No funding was received for conducting this study.

\medskip
\noindent\textbf{Competing interests.}
The author declares that he has no relevant
financial or non-financial interests to disclose.

\medskip
\noindent\textbf{Data availability.}
No datasets were generated or analyzed during
this study.

\medskip
\noindent\textbf{AI assistance.}
OpenAI's ChatGPT and Codex assisted with exploratory calculations,
checking proof steps, checking references, drafting and revising
the manuscript. The author takes responsibility for the final content.

\bibliographystyle{plainnat}
\bibliography{references}

@article{FinslerHadwiger,
  author = {Finsler, Paul and Hadwiger, Hugo},
  title = {Einige {Relationen} im {Dreieck}},
  journal = {Commentarii Mathematici Helvetici},
  volume = {10},
  pages = {316--326},
  year = {1937},
  doi = {10.1007/BF01214300}
}

@article{GiugiucDaoAltintas,
  author = {Giugiuc, Leonard Mihai and Dao, Thanh Oai and Altintas, Kadir},
  title = {An inequality related to the lengths and area of a convex quadrilateral},
  journal = {International Journal of Geometry},
  volume = {7},
  number = {1},
  pages = {81--86},
  year = {2018},
  url = {https://ijgeometry.com/wp-content/uploads/2018/04/81-86.pdf}
}

@misc{DaoQuestion,
  author = {Dao, Thanh Oai},
  title = {Strengthened version of isoperimetric inequality with {$n$}-polygon},
  howpublished = {MathOverflow, question 299056},
  year = {2018},
  note = {Posted April 30, 2018; includes an answer by Iosif Pinelis dated
          July 10, 2018. Accessed October 28, 2026},
  url = {https://mathoverflow.net/questions/299056}
}

@article{KuKu,
  author = {Ku, Hsu-Tung and Ku, Mei-Chin},
  title = {Analytic isoperimetric inequalities},
  journal = {Mathematical Inequalities \& Applications},
  volume = {3},
  number = {4},
  pages = {459--472},
  year = {2000},
  url = {https://files.ele-math.com/articles/mia-03-45.pdf}
}

@article{HittZhang,
  author = {Hitt, L. Richard and Zhang, Xin-Min},
  title = {Dynamic geometry of polygons},
  journal = {Elemente der Mathematik},
  volume = {56},
  number = {1},
  pages = {21--37},
  year = {2001},
  doi = {10.1007/S000170050086}
}

@article{IndreiNurbekyan,
  author = {Indrei, Emanuel and Nurbekyan, Levon},
  title = {On the stability of the polygonal isoperimetric inequality},
  journal = {Advances in Mathematics},
  volume = {276},
  pages = {62--86},
  year = {2015},
  doi = {10.1016/j.aim.2015.02.013}
}

@article{FisherRuoffShilleto,
  author = {Fisher, J. Chris and Ruoff, D. and Shilleto, J.},
  title = {Perpendicular polygons},
  journal = {The American Mathematical Monthly},
  volume = {92},
  number = {1},
  pages = {23--37},
  year = {1985},
  doi = {10.1080/00029890.1985.11971529}
}

@article{Indrei,
  author = {Indrei, Emanuel},
  title = {A sharp lower bound on the polygonal isoperimetric deficit},
  journal = {Proceedings of the American Mathematical Society},
  volume = {144},
  number = {7},
  pages = {3115--3122},
  year = {2016},
  doi = {10.1090/proc/12947}
}

@incollection{KourimskaSkuppinSpringborn,
  author = {Kou{\v r}imsk{\'a}, Hana and Skuppin, Lara and Springborn, Boris},
  title = {A variational principle for cyclic polygons with prescribed edge lengths},
  booktitle = {Advances in Discrete Differential Geometry},
  editor = {Bobenko, Alexander I.},
  publisher = {Springer},
  address = {Berlin, Heidelberg},
  pages = {177--195},
  year = {2016},
  doi = {10.1007/978-3-662-50447-5_5}
}

@article{Tang,
  author = {Tang, Dingyi},
  title = {Discrete {Wirtinger} and isoperimetric type inequalities},
  journal = {Bulletin of the Australian Mathematical Society},
  volume = {43},
  number = {3},
  pages = {467--474},
  year = {1991},
  doi = {10.1017/S0004972700029312}
}

@incollection{DemaineGassendORourkeToussaint,
  author = {Demaine, Erik D. and Gassend, Blaise and O'Rourke, Joseph and Toussaint, Godfried T.},
  title = {All polygons flip finitely\ldots\ right?},
  booktitle = {Surveys on Discrete and Computational Geometry: Twenty Years Later},
  editor = {Goodman, J. and Pach, J. and Pollack, R.},
  series = {Contemporary Mathematics},
  volume = {453},
  pages = {231--255},
  publisher = {American Mathematical Society},
  year = {2008},
  url = {https://erikdemaine.org/papers/Flips_DCG20/}
}

@article{Bogosel,
  author = {Bogosel, Beniamin},
  title = {Optimal {Finsler--Hadwiger} inequalities},
  journal = {Results in Mathematics},
  volume = {80},
  number = {3},
  pages = {Article 88},
  year = {2025},
  doi = {10.1007/s00025-025-02405-6}
}
\end{document}